\documentclass[pamm,a4paper,fleqn]{w-art}
\usepackage{times,cite,w-thm}
\usepackage[T1]{fontenc}
\usepackage[utf8]{inputenc}

\usepackage{graphicx, subcaption}
\usepackage{amsmath}
\usepackage{amsfonts}
\usepackage{tikz}
\newcommand*\diff{\mathop{}\!\mathrm{d}}
\newcommand\norm[1]{\left\lVert#1\right\rVert}

\begin{document}

\TitleLanguage[EN]
\title[Topology Optimization for Uniform Flow Distribution in Electrolysis Cells]{Topology Optimization for Uniform Flow Distribution in Electrolysis Cells}

\author{\firstname{Leon} \lastname{Baeck}\inst{1,}%
\footnote{Corresponding author: e-mail \ElectronicMail{leon.baeck@itwm.fraunhofer.de}}} 
\address[\inst{1}]{\CountryCode[DE]Fraunhofer Institute for Industrial Mathematics ITWM, Kaiserslautern, Germany}

\author{\firstname{Sebastian} \lastname{Blauth}\inst{1}}

\author{\firstname{Christian} \lastname{Leith\"auser}\inst{1}}

\author{\firstname{Ren\'e} \lastname{Pinnau}\inst{2}}
\address[\inst{2}]{\CountryCode[DE]RPTU Kaiserslautern-Landau, Technomathematics Group, Kaiserslautern, Germany}

\author{\firstname{Kevin} \lastname{Sturm}\inst{3}}
\address[\inst{3}]{\CountryCode[AT]TU Wien, Institute of Analysis and Scientific Computing, Vienna, Austria}

\AbstractLanguage[EN]
\begin{abstract}
In this paper we consider the topology optimization for a bipolar plate of a hydrogen electrolysis cell. We present a model for the bipolar plate using the Stokes equation with an additional drag term, which models the influence of fluid and solid regions. Furthermore, we derive a criterion for a uniform flow distribution in the bipolar plate. To obtain shapes that are well-manufacturable, we introduce a novel smoothing technique for the fluid velocity. Finally, we present some numerical results and investigate the influence of the smoothing on the obtained shapes.
\end{abstract}

\maketitle  

\section{Introduction}
Hydrogen plays an important role for achieving a climate neutral industry as well as climate neutral means of transportation. For that the hydrogen needs to be produced, for example, by means of hydrogen electrolysis cells. Typically, those cells split water into oxygen and the desired hydrogen by using (green) electrical energy. We focus on a special kind of electrolysis cell here, the so called proton exchange membrane (PEM) electrolysis cell \cite{Metz2023}, which works as follows: First, water enters the cell on the anode side and is distributed over the cell by the bipolar plate. The water wanders through the so-called porous transport layer (PTL) to the membrane where the electrochemical reaction takes place. The oxygen is then discharged from the cell again by the bipolar plate. Whereas positive hydrogen ions can pass through the membrane and form hydrogen at the cathode side. To guarantee the efficiency of such a PEM electrolysis cell, the flow through the anode side bipolar plate is essential. The water needs to be distributed uniformly over the whole cell so that the entire membrane area is reached. For more details on PEM electrolysis cells, we refer the reader, e.g., to \cite{Metz2023}.

Here, we only investigate the anode side of the bipolar plate. For simplicity, we neglect the electrochemical reaction and phase changes and only consider the distribution of a homogeneous fluid. To obtain a uniform flow through the bipolar plate, we use techniques from topology optimization based on topological sensitivity analysis \cite{Sokolowski1999Topology}.

Topology optimization considers the optimization of a domain by changing its topological features by either adding or removing material. It was first introduced in the context of solid mechanics, but has since been applied to a wide variety of applications, for example, in elasticity \cite{Eschenauer1994Topology, Amstutz2006new} and fluid mechanics \cite{Borrvall2003Topology, NSa2016Topological}. We use the concept of the topological derivative, which was first introduced in \cite{Eschenauer1994Topology} and later mathematically justified in \cite{Sokolowski1999Topology}. The topological derivative measures the sensitivity of a shape functional under infinitesimal topological changes. There are different approaches for solving topology optimization problems numerically and here we utilize a level-set approach which was introduced in \cite{Amstutz2006new}.

This paper is structured as follows. In Section 2, we introduce a model for the anode side bipolar plate using the Borvall-Petersson model from \cite{Borrvall2003Topology}. Additionally, we derive a criterion for a uniform flow distribution in the bipolar plate. To obtain geometries suitable for manufacturing, we consider a smoothing of the fluid velocity which is incorporated in the cost functional for the optimization. In Section 3, we summarize some basic concepts of topology optimization including the topological derivative as well as a gradient-based algorithm from \cite{Amstutz2006new}. Further, we formally present the topological derivative of the cost functional for our problem. Finally, in Chapter 4, we numerically solve the topology optimization problem and present our numerical results. Particularly, we focus on the influence of the smoothing on the optimized geometries. Our results show that topology optimization can be used to create novel bipolar plate designs.

\section{Model problem for the anode side bipolar plate}

Let us begin with presenting our mathematical model for the anode side of the bipolar plate as well as the corresponding topology optimization problem for a uniform flow distribution.

\subsection{The Borvall-Petersson model}

First, let the hold-all domain $D\subset\mathbb{R}^d$ be an open and bounded set. In this paper, we focus on the case $d=2$, but an extension to $d=3$ is also possible. Additionally, let $\Omega\subset D$ be measurable and open. Then, we identify $\Omega$ as the fluid region and $D\setminus\bar{\Omega}$ consequently as the solid part of $D$. To model the fluid flow, we use the Borvall Petersson model, which was first introduced in \cite{Borrvall2003Topology} for topology optimization of fluids. We assume that the boundary $\Gamma=\partial D$ of the hold-all domain $D$ is divided into three parts: The inlet $\Gamma_{\text{in}}$, where the water enters the plate with velocity $u_{\text{in}}$, the boundary $\Gamma_{\text{wall}}$, where we have a no-slip condition and $\Gamma_{\text{out}}$, where a natural outflow condition is applied. The model reads
\begin{equation}
	\begin{split}
	-\Delta u+ \alpha u + \nabla p &= 0 && \text{in } D,\\
	\text{div}(u) &= 0 && \text{in } D, \\
	u &= u_{\text{in}} && \text{on } \Gamma_{\text{in}}, \\
	u &= 0 && \text{on } \Gamma_{\text{wall}}, \\
	\partial_n u- pn &= 0 && \text{on } \Gamma_{\text{out}},
	\end{split}
	\label{stokesdarcy}
\end{equation}
where $u$ and $p$ denote the fluid velocity and pressure, respectively. Here, $n$ denotes the outer unit normal vector on $\Gamma$. To distinguish between the solid and fluid part of the domain the inverse permeability $\alpha$ is used, which is given by
\begin{equation*}
	\alpha(x) = \begin{cases} 
	\alpha_U &\text{if $x\in D\setminus\bar{\Omega}$},\\
	\alpha_L &\text{if $x\in\Omega$},
	\end{cases}
\end{equation*}
where $\alpha_L$ and $\alpha_U$ are positive constants. In particular, $\alpha$ is chosen small inside the fluid region and large in the solid part. For a detailed derivation of the model we refer to \cite{Borrvall2003Topology}.

\subsection{A criterion for uniform flow distribution}

To increase the efficiency of the electrolysis cell, we want to assure that the water gets distributed uniformly all over the cell. Without knowledge of the shape of $\Omega$, it is not obvious how to characterize a flow as uniform mathematically. Therefore, we take a different approach and introduce a threshold velocity magnitude $u_t > 0$. The goal is now to achieve a flow whose magnitude reaches at least the threshold. This criterion ensures that each part of the bipolar plate receives a sufficiently large flow and no dead spots occur. Additionally, as our model also yields, albeit small, flow velocities in the solid part, we have to restrict our criterion to the fluid part $\Omega$. The optimization goal then reads
\begin{equation}
	\label{optgoal1}
	\norm{u(x)}\geq u_t \text{  for all }x\in\Omega,
\end{equation}
where $\norm{\cdot}$ denotes the Euclidean norm on $\mathbb{R}^d$. 

There are two aspects to consider here. First, the target velocity constraint might not be fulfillable close to the boundary $\partial\Omega$ of the fluid domain due to small velocities in the solid part. Additionally, as we do not have any restrictions on the solid part, large solid areas might appear in the design structures. In fact, we observed this behavior in numerical tests. Usually, we want to avoid solid structures above a certain size as the flow inside the porous transport layer is hindered by such objects and this leads to a degradation of the cell efficiency.

\subsection{Smoothing of the velocity}

To avoid large solid structures and diminish the effect of small flow velocities in the solid parts, our goal is to extend the target velocity goal $(\ref{optgoal1})$ to the hold-all domain $D$. We consider a smoothed flow velocity which is defined on the entire hold-all domain D. For that we use the smoothing characteristics of the heat equation. We consider the following equation
\begin{equation}
	\begin{split}
	\partial_tu_s-\Delta u_s &= 0 && \text{in } D,\\
	\nabla u_s\cdot n &= 0 && \text{on } \partial D, \\
	u_s(x,0) &= u(x) && \text{in } D.
	\end{split}
	\label{heat}
\end{equation}
Here, we have the velocity solution $u$ of the Stokes Darcy equation $(\ref{stokesdarcy})$ as the initial condition. To reduce the numerical effort, we discretize $(\ref{heat})$ with a single implicit Euler step with step length $\Delta t > 0$ and arrive at 
\begin{equation}
	\begin{split}
	\frac{u_s-u}{\Delta t}-\Delta u_s &= 0 && \text{in } D,\\
	\nabla u_s\cdot n &= 0 && \text{on } \partial D.
	\end{split}
	\label{heat1step}
\end{equation}
By modifying the time step length $\Delta t$ we can control the influence of the smoothing, where the smoothing effect gets larger for increasing step lengths. With the smoothed velocity $u_s$, we can now extend the criterion $(\ref{optgoal1})$ to the entire hold-all domain by considering
\begin{equation}
	\label{optgoal2}
	\norm{u_s(x)}\geq u_t \text{  for all }x\in D.
\end{equation}
Finally, we apply the Moreau-Yosida regularization \cite{Hinze2009Optimization} to $(\ref{optgoal2})$ to arrive at our objective function
\begin{equation}
\label{objective}
J(\Omega,u) = \int_{D} \min \left( 0, \norm{u_s} - u_d \right)^2 \diff x.
\end{equation}
\begin{remark}
	Using the smoothed velocity has several advantages for the topology optimization. First, when we consider a small solid inclusion, we expect that it does not severely effect the flow in the underlying porous transport layer. In the smoothed velocity $u_s$ with a suitable time step $\Delta t$, the effect of the small solid inclusion will not be resolved anymore, so that the smoothed velocity $u_s$ is not influenced anymore by the small velocity inside the solid region and the target velocity goal $(\ref{optgoal2})$ will be fulfilled. On the other hand, if we consider a large solid inclusion, it will still be resolved by the smoothed velocity $u_s$ and the criterion $(\ref{optgoal2})$ will not be satisfied. Therefore, large solid structures can be avoided with this approach.
\end{remark}

\subsection{Topology optimization problem}

We summarize $(\ref{stokesdarcy})$, $(\ref{heat1step})$ and $(\ref{objective})$ to state our topology optimization problem
\begin{equation}
	\begin{split}
	\min_{\Omega}J(\Omega,u)&= \int_{D} \min \left( 0, \norm{u_s} - u_t \right)^2 \diff x \\
	&\text{such that } (\ref{stokesdarcy}) \text{ as well as } (\ref{heat1step}) \text{ are fulfilled in the weak sense}, \\
	&V_L \leq |\Omega| \leq V_U.
	\end{split}
	\label{optproblem}
\end{equation}
Here, we introduce an additional volume constraint for the fluid volume with $V_L$ and $V_U$ being the lower and the upper border, respectively. The goal is now to perform a topology optimization for problem $(\ref{optproblem})$ using the topological derivative.

\section{Topology optimization}

We seek to apply a gradient-based topology optimization approach for $(\ref{optproblem})$. For that we introduce the topological derivative and state the corresponding topological derivative for the objective $(\ref{objective})$. Furthermore, we briefly present the gradient-based topology optimization algorithm we use for our numerical experiments.

\subsection{Topological derivative}

We consider an analogous setting to before, where $D$ is an open subset of $\mathbb{R}^d$. We define $\mathcal{A}=\{\Omega\subset D \text{ s.t. } \Omega \text{ open}\}$ as the set of all open subsets of $D$. To present the topological derivative, we consider a general shape functional $S$, given by
\begin{equation*}
	\Omega\in\mathcal{A}\mapsto S(\Omega)\in\mathbb{R}.
\end{equation*}
Additionally, we assume to have a fixed shape $\Omega$ and a fixed point $z\in D\setminus\partial\Omega$. The topological derivative then measures the sensitivity of the shape functional $S$ with respect to infinitesimal topological changes of the shape $\Omega$ around the point $z$. We characterize the perturbed shape by
\begin{equation*}
	\Omega_{z,\epsilon}=\begin{cases}
	\Omega\setminus\bar{\omega}_{z,\epsilon}, & z\in\Omega\\
	\Omega\cup\omega_{z,\epsilon}, & z\in D\setminus\bar{\Omega},
	\end{cases}
\end{equation*}
where $\omega_{z,\epsilon}=z+\epsilon\omega$ with $\omega\subset\mathbb{R}^d$ and $0\in\omega$ represents the shape of the perturbation. Furthermore, we assume to have a positive function $l$ with $\lim_{\epsilon\rightarrow0}l(\epsilon)=0$. The topological derivative, see for example \cite{Sokolowski1999Topology}, is then defined as
\begin{equation}
	D_TS(\Omega,\omega)(z)=\lim_{\epsilon\rightarrow0}\frac{S(\Omega_{z,\epsilon})-S(\Omega)}{l(\epsilon)},
\end{equation}
if the limit exists. Typical choices are $w=B_1(0)$ with $l(\epsilon)=|\omega_{z,\epsilon}|$, where $B_1(0)$ is the open unit ball with center $0$.

\subsection{Topological derivative for our model problem}
We formally present the topological derivative of our model problem $(\ref{optproblem})$, which can be derived, for example, with an averaged adjoint approach [8]. The topological derivative for $(\ref{objective})$ is given by
\begin{equation}
	\label{topder}
	D_TJ(z)=-(\alpha_U-\alpha_L)u(z)v(z)
\end{equation}
for all $z\in D\setminus\partial\Omega$. Here, $u$ is the weak solution of $(\ref{stokesdarcy})$ and $v$ is the adjoint velocity which solves
\begin{equation}
\begin{split}
-\Delta v+ \alpha v + \nabla q-\frac{1}{\Delta t}v_s &= 0 && \text{in } D,\\
\text{div}(q) &= 0 && \text{in } D, \\
v &= 0 && \text{on } \Gamma_{\text{in}}, \\
v &= 0 && \text{on } \Gamma_{\text{wall}}, \\
\partial_n v- qn &= 0 && \text{on } \Gamma_{\text{out}},
\end{split}
\label{stokesadj}
\end{equation}
and $v_s$ is the adjoint of the smoothed velocity which solves the equation
\begin{equation}
	\label{heatadjoint}
	\begin{split}
	\frac{1}{\Delta t}v_s-\Delta v_s &= 2\frac{u_s}{\norm{u_s}}\min(0,\norm{u_s}-u_t) && \text{in }D,\\
	\nabla v_s\cdot n &= 0 &&\text{on }\partial D.
	\end{split}
\end{equation}
The actual computation of the topological derivative is beyond the scope of this paper and a topic of future research.  For more details regarding topological sensitivity analysis, we refer the reader to \cite{Sokolowski1999Topology}.

\subsection{A gradient-based solution algorithm for topology optimization}

We present an algorithm using the topological derivative as a search direction in order to solve topology optimization problems. The method was introduced by Amstutz and Andr\"a \cite{Amstutz2006new}. For the sake of brevity, we only give a short presentation here and refer the reader to \cite{Blauth2023Quasi} for an overview over established topology optimization algorithms as well as novel quasi-Newton methods for topology optimization. First, we introduce a continuous level-set function $\psi:D\rightarrow\mathbb{R}$ to characterize the fluid as well as the solid part of the hold-all domain $D$ by
\begin{equation*}
	\psi(z)\begin{cases}
	<0, & z\in\Omega,\\
	=0, & z\in\partial\Omega,\\
	>0. & z\in D\setminus\bar{\Omega}.
	\end{cases}
\end{equation*}
Additionally, we define the generalized topological derivative
\begin{equation*}
g(z)=\begin{cases}
-D_TS(\Omega,\omega)(z), &z\in\Omega,\\
+D_TS(\Omega,\omega)(z), &z\in D\setminus\bar{\Omega}.
\end{cases} 
\end{equation*}
The algorithm then updates the level-set function consequently by a linear combination of itself and the generalized topological derivative on the $L^2$-sphere. The update formula reads
\begin{equation}
\label{update}
\psi_{n+1}=\frac{1}{\sin(\theta_n)}\left[\sin((1-\theta_n)\kappa_n)\psi_n+\sin(\theta_n\kappa_n)\frac{g_n}{\norm{g_n}_{L^2(D)}}\right].
\end{equation}
Here, $\theta_n$ is the $L^2$-angle between the level-set function and the generalized topological derivative
\begin{equation*}
\theta_n = \arccos\left[\frac{\langle g,\psi\rangle_{L^2(D)}}{\norm{g}_{L^2(D)} \norm{\psi}_{L^2(D)}}\right].
\end{equation*}
Additionally, $\kappa_n$ is chosen by a line search approach to guarantee that the objective function value decreases each iteration. This procedure is carried out until a local optimality condition is reached. Such an optimality condition has been derived e.g. in \cite{Amstutz2011Analysis}. Numerically, the algorithm is stopped, when the angle $\theta_n$ becomes smaller than a certain numerical tolerance $\epsilon_\theta$. This then implies a local optimality condition, again the reader is referred to \cite{Amstutz2006new} or \cite{Amstutz2011Analysis}. We want to conclude by stating that all $\psi_n$ fulfill $\norm{\psi_n}_{L^2(D)}=1$ if the condition holds for the initial iteration.

\section{Numerical results}

We turn back to optimization problem $(\ref{optproblem})$ and describe our numerical setting. First, the hold all domain is given by $D=(0,1)\times(0,1)$. We define a quadratic inflow profile $u_{\text{in}}$ given by 
\begin{equation*}
	u_{\text{in}}=\begin{bmatrix}
	-\frac{400}{9}(y-0.35)(y-0.65)\\
	0
	\end{bmatrix} \; \text{for }x=0 \text{ and } 0.35\leq y\leq 0.65
\end{equation*}
on the inflow boundary $\Gamma_{\text{in}}$. The situation is displayed in Figure \ref{fig_model}. The inverse permeability $\alpha$ for the fluid as well as the solid part as is defined as
\begin{equation*}
\alpha_L=\frac{2.5}{125^2}, \;\; \alpha_U=\frac{2.5}{0.008^2}.
\end{equation*}
Furthermore, the stopping criterion for the algorithm is given by $\epsilon_\theta=0.035$ which is equivalent to an angle of $2$ degrees. The target velocity is chosen as $u_t=0.1$ and the lower and upper values for the volume constraint are given by $V_L=0.5$ and $V_U=0.7$, respectively.

\begin{vchfigure}
	\begin{tikzpicture}[scale=0.3]
	\draw[line width = 0.35mm] (0,0) -- (10,0) -- (10,10) -- (0,10) -- (0,0);
	\draw[black, line width = 0.30mm]   plot[smooth,domain=3.5:6.5] ({-(\x-3.5)*(6.5-\x)}, \x);
	\draw[black, line width = 0.30mm]   plot[smooth,domain=3.5:6.5] ({(\x-3.5)*(6.5-\x)+10}, \x);
	\draw[black, line width = 0.30mm, ->] (-1.6875,5.75) -- (0,5.75);
	\draw[black, line width = 0.30mm, ->] (-2.25,5) -- (0,5);
	\draw[black, line width = 0.30mm, ->] (-1.6875,4.25) -- (0,4.25);
	\draw[black, line width = 0.30mm, ->] (10,5.75) -- (11.6875,5.75);
	\draw[black, line width = 0.30mm, ->] (10,5) -- (12.25,5);
	\draw[black, line width = 0.30mm, ->] (10,4.25) -- (11.6875,4.25);
	\draw[black, line width = 0.30mm, <->] (0.5,5) -- (0.5,10);
	\draw[black, line width = 0.30mm, <->] (0.5,0) -- (0.5,3.5);
	\draw[black, line width = 0.30mm, <->] (0,10.5) -- (10,10.5);
	\draw[black, line width = 0.30mm, <->] (13,0) -- (13,10);
	\node[] (a) at (5,5) {$D$};
	\node[] (b) at (15.25,5) {$l=1.0$};
	\node[] (b) at (5,11.25) {$l=1.0$};
	\node[] (b) at (1.75,1.75) {$0.35$};
	\node[] (b) at (1.5,7.5) {$0.5$};
	\end{tikzpicture}
	\captionsetup{margin=4cc}
	\caption{Schematic setup for the bipolar plate.}
	\label{fig_model}
\end{vchfigure}

\subsection{Numerical implementation}

We give a short overview of the software used for the numerical implementation. First, all underlying PDEs are solved using the finite element package FEniCS \cite{Alnes2015FEniCS}. We use a uniform triangular mesh consisting of 10.000 elements. The Stokes Darcy equation $(\ref{stokesdarcy})$ is discretized using LBB-stable Taylor-Hood finite elements which consist of quadratic Lagrange elements for the velocity and linear Lagrange elements for the pressure. Additionally, for the smoothing equation $(\ref{heat1step})$, continuous quadratic Lagrange elements are used.

Furthermore, we use the software package cashocs \cite{Blauth2023Cashocs} for adjoint computations. Cashocs is an open source software package with a python interface which is based on FEniCS. The software can be used to solve optimization problems constrained by partial differential equations in the context of shape optimization and optimal control in an automated fashion. Furthermore, cashocs implements a discretization of the continuous adjoint approach which allows for adjoints to be computed automatically.

As we also introduced a volume constraint for $(\ref{optproblem})$, we give an explanation on how the constraint is handled numerically. After each update of the level-set function we perform a projection of the level-set function onto the set of the admissible shapes. First, we assume to have a newly computed update of the level-set function $\psi_{n+1}$, which was generated with $(\ref{update})$, and the corresponding shape $\Omega_{n+1}$. Additionally, we assume that the volume constraint is not fulfilled, otherwise nothing needs to be done. Without loss of generality we assume that $|\Omega_{n+1}|>V_U$. To perform the projection, we move the level-set function upwards by a strictly positive constant until the upper border $V_U$ of the volume constraint is reached. This means we search for a $c>0$, to move the level-set function $\tilde{\psi}_{n+1}=\psi_{n+1}+c$, such that the corresponding shape $\tilde{\Omega}_{n+1}$ fulfills $|\tilde{\Omega}_{n+1}|=V_U$. The case $|\Omega_{n+1}|<V_L$ is handled analogously by taking a strictly negative constant to alter the level-set function. Numerically, we compute the constant $c$ by a bisection approach with numerical accuracy $\epsilon_c=10^{-4}$. As we do not need to solve any additional PDEs, the computational cost for the projection procedure is reasonable.

\subsection{Numerical results}

\begin{figure}
	\includegraphics[width=.3\linewidth]{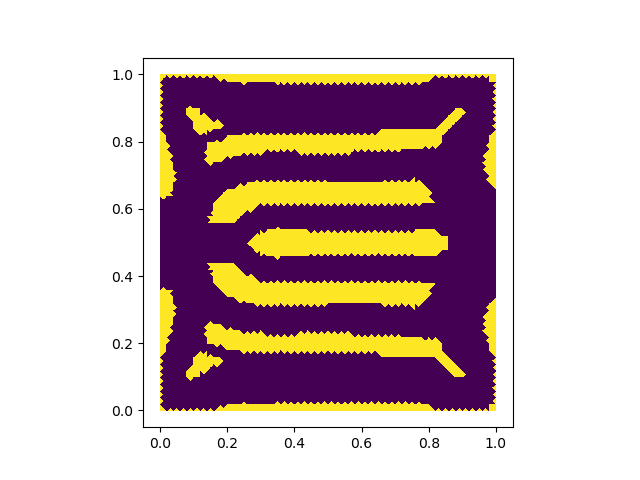}~a)
	\hfil
	\includegraphics[width=.3\linewidth]{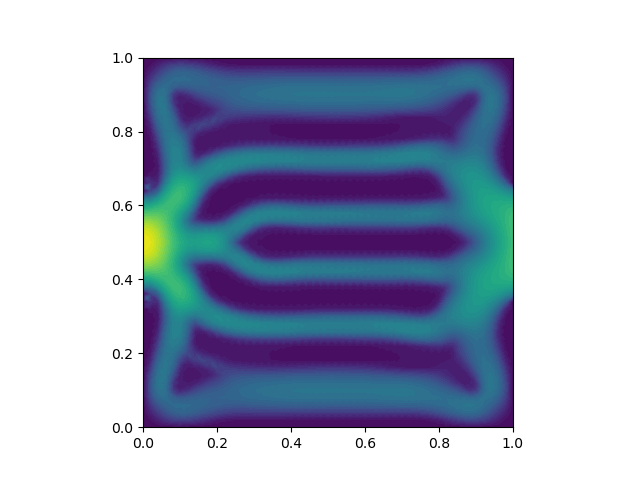}~b)
	\hfil
	\includegraphics[width=.3\linewidth]{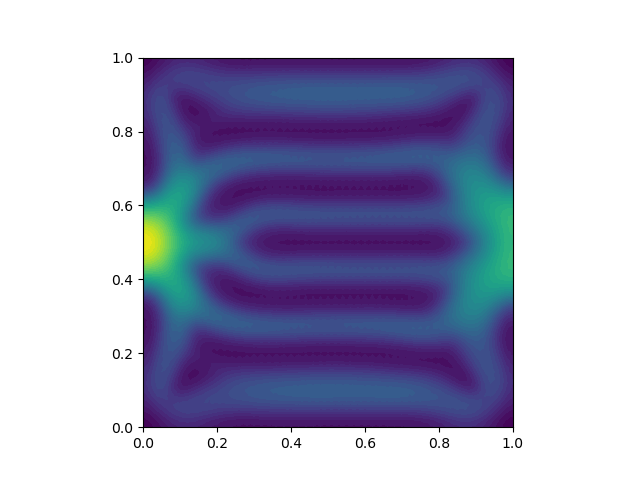}~c)
	\caption{Result of the topology optimization problem $(\ref{optproblem})$ with $\Delta t=10^{-3}$. In Figure \textbf{a} the final shape is displayed, in Figure \textbf{b} the corresponding norm of the velocity field and in Figure \textbf{c} the norm of the smoothed velocity field.}
	\label{fig:result1}
\end{figure}

We start by choosing $\Delta t=10^{-3}$ as the step length for the smoothing equation $(\ref{heat1step})$. The resulting optimized geometry and flow velocity are displayed in Figure \ref{fig:result1}. The algorithm took 137 iterations to reach the target accuracy. The objective function value is given by $8.6\cdot 10^{-8}$. From the obtained velocity fields and the previous investigations, we observe that the optimization algorithm indeed produced a geometry with a very uniform flow distribution. To do so, six fluid channels have been formed and the flow is distributed rather uniformly between them. In more detail, the target velocity goal $(\ref{optgoal2})$ for the smoothed velocity $u_s$ is reached on $93.84\%$ of the whole domain. Furthermore, this extends to a $91.79\%$ fulfillment of $(\ref{optgoal1})$ for $u$ on the fluid area.  Additionally, all obstacle sizes in fluid flow direction are about the same order of magnitude. 

\begin{figure}
	\includegraphics[width=.3\linewidth]{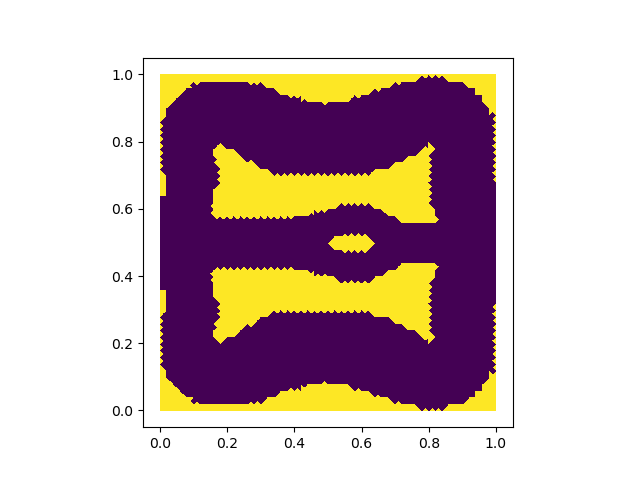}~a)
	\hfil
	\includegraphics[width=.3\linewidth]{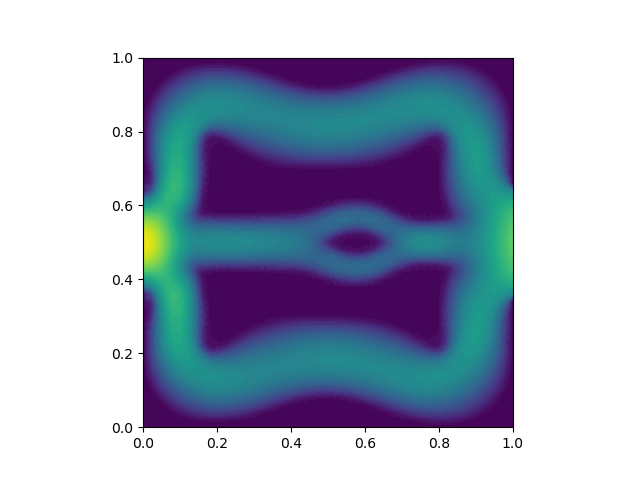}~b)
	\hfil
	\includegraphics[width=.3\linewidth]{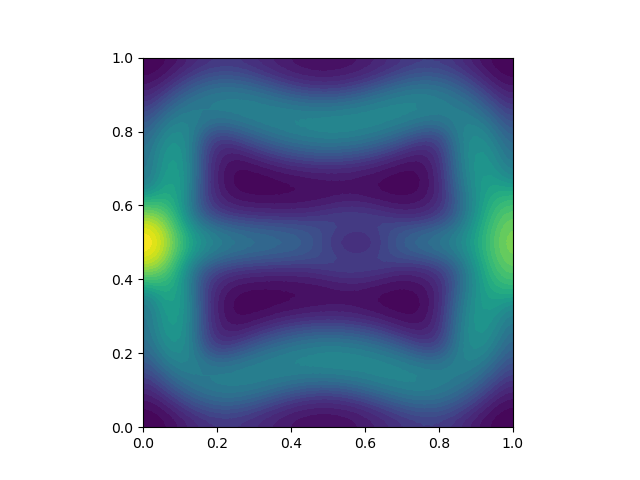}~c)
	\caption{Result of the topology optimization problem $(\ref{optproblem})$ with $\Delta t_1=5\cdot 10^{-3}$. In Figure \textbf{a} the final shape is displayed, in Figure \textbf{b} the corresponding norm of the velocity field and in Figure \textbf{c} the norm of the smoothed velocity field.}
	\label{fig:result2}
\end{figure}

\begin{figure}
	\includegraphics[width=.3\linewidth]{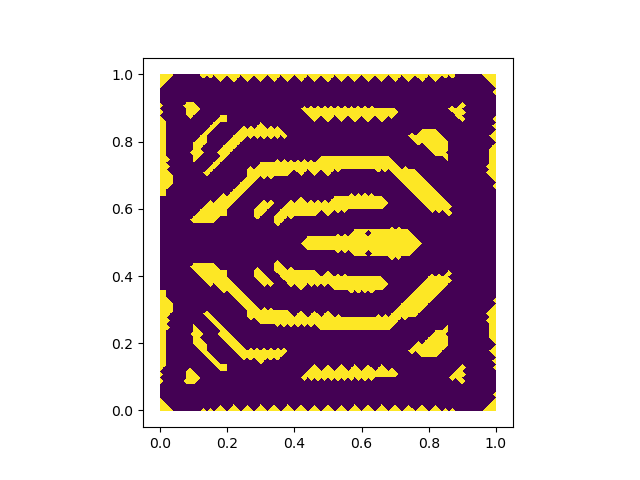}~a)
	\hfil
	\includegraphics[width=.3\linewidth]{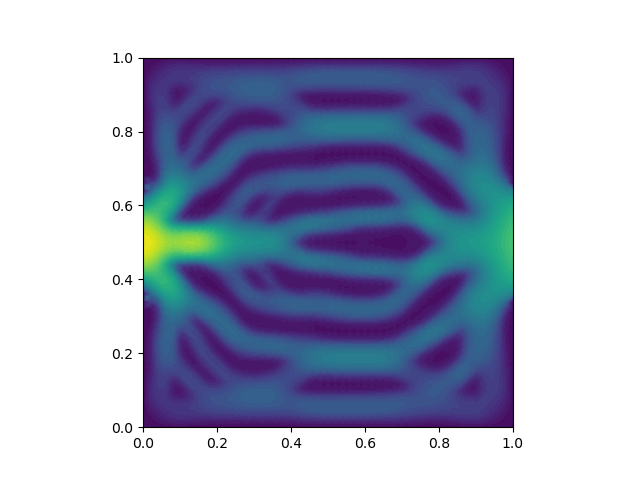}~b)
	\hfil
	\includegraphics[width=.3\linewidth]{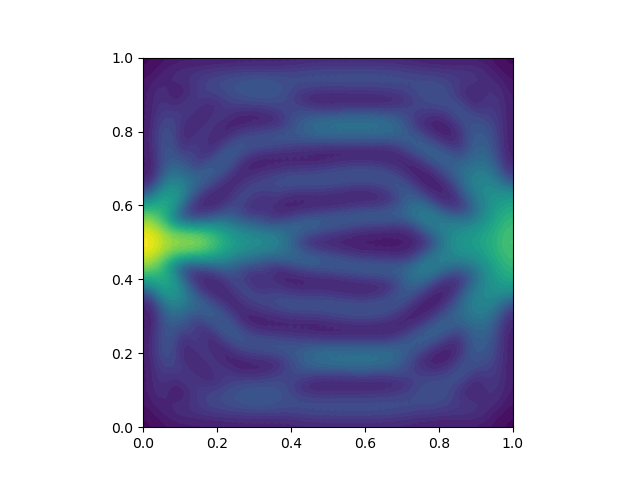}~c)
	\caption{Result of the topology optimization problem $(\ref{optproblem})$ with $\Delta t_2=5\cdot 10^{-4}$. In Figure \textbf{a} the final shape is displayed, in Figure \textbf{b} the corresponding norm of the velocity field and in Figure \textbf{c} the norm of the smoothed velocity field.}
	\label{fig:result3}
\end{figure}

We want to investigate the obstacle sizes for different values of the smoothing parameter $\Delta t$. For that we choose two additional values for the step length given by $\Delta t_1=5\cdot 10^{-3}$ as well as $\Delta t_2=5\cdot 10^{-4}$. The results can be seen in Figure \ref{fig:result2} for $\Delta t_1$ and in Figure \ref{fig:result3} for $\Delta t_2$. It is apparent that the solid structures for the higher time step length get larger which can be explained by the higher influence of the smoothing. The larger the smoothing parameter is chosen, the larger a solid inclusion has to be in order to be resolved in the smoothed velocity $u_s$. In particular, we now observe that only 3 major fluid channels can be found in the optimized geometry. Here, the target velocity goal $(\ref{optgoal2})$ is fulfilled on $98.26\%$ of the domain $D$. For the smaller time step length $\Delta t_2$ on the contrary, the solid structures get smaller, which is again a consequence of the smaller influence of the smoothing and we observe that the geometry branches out into several channels which are then merged back together towards the outlet. We achieve a $89.46\%$ realization of the target velocity goal for the smoothed velocity $u_s$.

\section{Conclusion}

Summarizing, we presented a model for reaching a uniform flow distribution in the anode side bipolar plate for hydrogen electrolysis cells. We showed in numerical tests that a uniform flow distribution is indeed reached in the optimized shapes. Additionally, we investigated the influence of the time step length $\Delta t$ on the size of the solid structures, giving us a tool to control the object sizes in the resulting shapes.


\begin{acknowledgement}
  L. Baeck gratefully acknowledges financial support from the Fraunhofer Institute of Industrial Mathematics ITWM.
\end{acknowledgement}

\vspace{\baselineskip}

\end{document}